\documentclass[a4paper]{easychair}
\usepackage{amssymb}
\newcommand{\player}[1]{\mathbf{#1}}
\newcommand{\Opp}[0]{\player{O}}
\newcommand{\Pro}[0]{\player{P}}
\newcommand{\bottom}[0]{\perp}
\newcommand{\system}[1]{\textsf{#1}}
\newcommand{\kuno}[0]{\system{Kuno}}

\newcommand{\tptpfourx}{\system{TPTP4X}}

\usepackage{xcolor}

\begin{document}
\title{Dialogues for proof search}
\titlerunning{Dialogues for proof search}
\author{Jesse Alama\inst{1}\thanks{Supported by FWF grant P25417-G15 (LOGFRADIG).}}
\institute{Theory and Logic Group{\\}Technical University of Vienna{\\}Vienna, Austria{\\}\email{alama@logic.at}}
\authorrunning{Alama}

\newtheorem{theorem}{Theorem}
\newtheorem{conjecture}{Conjecture}
\newtheorem{definition}{Definition}

\clearpage

\maketitle

\begin{abstract}Dialogue games are a two-player semantics for a variety of logics, including intuitionistic and classical logic.  Dialogues can be viewed as a kind of analytic calculus not unlike tableaux.  Can dialogue games be an effective foundation for proof search in intuitionistic logic (both first-order and propositional)?  We announce {\kuno}, an automated theorem prover for intuitionistic first-order logic based on dialogue games.\end{abstract}

\pagestyle{empty}

\section{Introduction}\label{sect:introduction}

Dialogue games (``dialogues'' for short) are a game-theoretic semantics for intuitionistic logic.  The game starts with a logical formula $\phi$ asserted by the Proponent ($\Pro$), who takes the stance that the formula is valid, against the Opponent ($\Opp$) who disputes this.  The players take turns and move according to certain rules.  As with other game-based semantics for logics, the focus is less at the level of a particular play of the game (sequence of moves that follow the rules) and more at the level of strategy (a way of playing the game to ensure a certain outcome).  The main result about dialogue games that is of interest for us here is:

\begin{theorem}\label{felscher-theorem}A formula $\phi$ is intuitionistically valid iff there exists a winning strategy for $\phi$.\end{theorem}

(A winning strategy for $\phi$ is a way of playing a dialogue game such that, for every move that $\Opp$ can make, $\Pro$ can respond in such a way that a win [for $\Pro$] is ensured.)

The dialogical approach to logic (``dialogical logic'', for short) is found mainly in philosophical discussions about logic and semantics.  Dialogues differ from other game-based approaches to logic, such as Hintikka-style games.  There, for instance, players have certain roles (e.g., ``Abelard'' and ``Eloise'') that can switch as the game proceeds; in dialogue games, the two players $\Pro$ and $\Opp$ play the same ``role'' throughout the game (they do not swap sides).

The question in focus in the present paper is: \emph{Can dialogue games be profitably understood as a proof search calculus?}  The aim is to cast dialogues in a new light, by implementing them and casting them in a light that they do not normally see.  We also hope to reap fresh insights into the foundations of dialogical logic by seeing how they respond to the pressures, so to speak, of everyday proof search.  The fruit of our efforts is {\kuno}, an automated theorem prover for first-order intuitionistic logic that is based on dialogue games.  Section~\ref{sec:rules} discusses the game calculus.  We shall see that dialogue games are essentially a kind of analytic calculus with similarities to tableaux.  Section~\ref{sec:implementation} briefly discusses the implementation and operation of {\kuno}.  Section~\ref{sec:evaluation} evaluates {\kuno} on a part of the Intuitionistic Logic Theorem Proving Library (ILTP), a collection of theorem proving problems~\cite{raths2005iltp,raths2007iltp}.  Section~\ref{sec:conclusion} concludes with further problems to be solved.

{\kuno} is available online at \url{http://github.com/jessealama/dialogues}.

\section{Dialogue games for intuitionistic first-order logic}\label{sec:rules}

For a thorough introduction to dialogical logic, see~\cite{sep-logic-dialogical} or~\cite{rahman2005how}.  We work in a first-order language built from $\forall$, $\exists$, $\neg$, $\vee$, $\wedge$, and $\rightarrow$.  It is assumed that equality is not present, nor $\bottom$, nor $\top$.  (All of these restriction can be lifted by suitable rewritings, but for the sake of simplicity we give the traditional formulation of dialogues.)  Dialogue games involve not just formulas, but also so-called \emph{symbolic attacks} $?$ (akin to asking ``which?''), $\wedge_{L}$, and $\wedge_{R}$ (akin to prompting the other player ``defend the left-hand side'' or ``defend the right-hand side'').  Together formulas and symbolic attacks are called \emph{statements}; they are what is asserted by $\Pro$ and $\Opp$ in a dialogue game.

The rules governing dialogues are divided into two types.  \emph{Particle rules} can be seen as specifying the meaning of connectives in a local fashion and say how formulas can be attacked and defended depending on their main connective.  By contrast, \emph{structural rules} operate globally and define what sequences of attacks and defenses count as dialogues, thus giving a kind of global meaning to the connectives.

\begin{table}[t]
\begin{center}
  \setlength{\tabcolsep}{5pt}
  \begin{tabular}{c|c|c}
    \textbf{Assertion} & \textbf{Attack} & \textbf{Response}\\
    \hline
    $\phi \wedge \psi$       & $\wedge_{L}$ & $\phi$\\
                             & $\wedge_{R}$ & $\psi$\\
    $ \phi \vee \psi$        & $?$         & $\phi$ or $\psi$\\
    $ \phi \rightarrow \psi$ & $\phi$      & $\psi$\\
    $ \neg\phi$              & $\phi$      & ---\\
    \hline
    $\forall x \phi$         & $?_{t}$     & $\phi[t/x]$\\
    $\exists x \phi$         & $?$        & $\phi[t/x]$\\
    \hline
  \end{tabular}
\medskip
  \caption{Particle rules for dialogue games}
  \label{tab:particle-rules}
\end{center}
\end{table}

Dialogue games start with the Proponent ($\Pro$) making the initial assertion.  Play alternates between $\Pro$ and the Opponent $\Opp$.  Every move is either an attack on something previously asserted or a defense against an attack. The standard particle rules are given in Table \ref{tab:particle-rules}. According to the first row, there are two possible attacks against a conjunction: The attacker specifies whether the left or the right conjunct is to be defended, and the defender then continues the game by asserting the specified conjunct.  The second row says that there is one attack against a disjunction; the defender then chooses which disjunct to assert.  The interpretation of the third row is straightforward.  The fourth row says that there is no way to defend against the attack against a negation; the only appropriate ``defense'' against an attack on a negation $\neg\phi$ is to continue the game with the new information $\phi$.  The particle rule for $\forall$ says that the challenger picks an instance (a term) and it is up to the original claimant to defend the instance of the universal generalization.  For $\exists$, the challenge is simply: ``which one?''  The way to proceed is to pick an instance of the existential generalization.

The structural rules are:

\begin{itemize}
\item $\Pro$ may assert an atomic formula only if $\Opp$ has asserted it earlier.
\item Only the most recent open attack may be defended.  (An attack is open if there there is no defense against it.)
\item Attacks may be defended only once.
\item $\Pro$'s assertions may be attacked at most once.
\end{itemize}

The game ends if no possible move can be made.  If $\Opp$ cannot move, then it is said that $\Pro$ has won the game; if $\Pro$ cannot move, then $\Opp$ has won the game.  (It is possible, even at the propositional level, for games to go on infinitely, with neither player winning.)

In addition to these standard structural rules, another rule is often considered in connection with dialogues:

\begin{description}
\item[E] $\Opp$ must immediately respond to $\Pro$'s moves.
\end{description}

For brevity, by ``the E rules'' we mean the structural rules together with the E rule.  The standard name in the dialogue game literature for the structural rules presented here is ``D''.  Evidently, when the so-called E rule is present, $\Opp$ is rather tightly constrained.  A consequence of the E rule being present is that whenever $\Pro$ defends against an attack, $\Opp$ must immediately attack $\Pro$'s move.

\begin{theorem}[Felscher]There exists a winning strategy for $\phi$ iff there exists a winning strategy for $\phi$ that adheres to the E rules.\end{theorem}

(Note that our assumption that E is present is helpful for proof search considerations.  It is another matter to philosophically justify the inclusion of E.  We are relying on the fact that dialogue validity is the same with or without E, a result proved by Felscher.)








\section{Implementation}\label{sec:implementation}

{\kuno} is a Common Lisp (CL) program.  One can run {\kuno} within a CL Read-Eval-Print loop (REPL), or from the commandline by first compiling the CL sources.  At the moment, the only tested CL implementation is SBCL (Steel Bank Common Lisp), a major open-source CL implementation (\url{http://www.sbcl.org}).

The name ``Kuno'' is a tribute to Kuno Lorenz, one of the foundational figures of the field of dialogue games~\cite{kuno2001basic}.

{\kuno} is based on a previous program that was designed to support a kind of interactive proof search, with a web-based frontend, is used~\cite{DBLP:conf/lpar/AlamaU10a}.

Given a TPTP problem, we first separate the conjecture formula from the non-conjecture formula.  (It is an error to invoke {\kuno} on a TPTP problem that lacks a single conjecture formula.)  If the TPTP input contains only a conjecture formula, then this evidently is how the game shall begin.  Otherwise, we form in the obvious way an implication whose consequent is the conjecture formula and whose antecedent is the conjunction of the non-conjecture formulas.

\section{Evaluation on ILTP}\label{sec:evaluation}

We consider first the propositional part of the ILTP (version 1.1.2, available at \url{http://www.cs.uni-potsdam.de/ti/iltp/formulae.html}).  Table~\ref{tab:e-evaluation} contains the result of working with propositional problems in the LCL (devoted to Logic Calculi) and SYN (Syntactic) sections of the ILTP library.  We developed strategies to a depth limit of 30 (that is, if a strategy ever exceeded depth 30, it was discarded from the search even if potentially it could be completed to a winning strategy).  ``Depth'' means that {\kuno} did terminate, but the best it can say is that there is no strategy below the given depth limit.  ``Timeout'' means that computation had to be halted by a time limit.  ``Crash'' means crash.


\begin{table}
\begin{tabular}{|c|c|c|l|}
Problem      & Intended SZS Status  & Computed SZS Status agrees? & Reason\\
\hline\hline
LCL181+1     & Non-Theorem & - & depth   \\
LCL230+1     & Non-Theorem & - & depth   \\
\hline
SYN001+1     & Non-Theorem & - & depth   \\
SYN007+1.014 & Non-Theorem & - & crash   \\
SYN040+1     & Non-Theorem & - & timeout \\
SYN041+1     & Theorem     & + &         \\
SYN044+1     & Theorem     & + &         \\
SYN045+1     & Theorem     & + &         \\
SYN046+1     & Non-Theorem & - & depth   \\
SYN047+1     & Non-Theorem & - & timeout \\
SYN387+1     & Non-Theorem & - & depth   \\
SYN388+1     & Non-Theorem & - & depth   \\
SYN389+1     & Non-Theorem & - & depth   \\
SYN390+1     & Theorem     & + &         \\
SYN391+1     & Theorem     & + &         \\
SYN392+1     & Non-Theorem & - & timeout \\
SYN393+1     & Non-Theorem & - & timeout \\
SYN416+1     & Non-Theorem & - & depth   \\
SYN915+1     & Theorem     & + &         \\
SYN916+1     & Non-Theorem & + &         \\
SYN977+1     & Non-Theorem & - & depth   \\
SYN978+1     & Theorem     & + &         \\
\hline
\end{tabular}
\caption{\label{tab:e-evaluation}An evaluation of the E ruleset on several problems from the ILTP library (propositional part)}
\end{table}

The initial experiment was useful not only for identifying bugs, but for gaining additional insight into dialogues as a decision procedure for propositional intuitionistic logic.  In the cases where a result is indeed a theorem, we were not especially surprised, in view of previous experience with the architecture underlying {\kuno}, which was focused on working with known theorems.  It was the case of non-theorems where we were more interested.  We found that an important obstacle that prevents the program from terminating with useful information is the possibility of endless repetition or duplication by one of the players.  In the case where we are dealing with a formula $\phi$ that is not intuitionistically valid, we find that there are two possibilities:

\begin{itemize}
\item $\Opp$ is able to repeat moves ad infinitum.
\item The dialogue search tree (a complete development of all possible dialogues whatsoever starting from an initial formula) is finite, but it contains no winning strategy.
\end{itemize}




The second case is generally detected by {\kuno}; the first case is more interesting.  Motivated by such concerns, we are led to consider an additional ``no repetition'' constraint.

\begin{description}
\item[No-Repeats] Neither player can repeat a move.
\end{description}

We are at the moment unable to prove this constraint preserves completeness in the context of intuitionistic logic (though an analogous restriction preserves completeness in the case of classical logic~\cite{clerbout2013firstorder}).





Happily, when working with the E ruleset with the No-Repeats constraint, we are able to solve all of the unsolved problems of Table~\ref{tab:e-evaluation} except SYN007+1.014, SYN047+1, and SYN393+1.  Looking into more detail on the failure to come a decision on these three problems, we find two sources of difficulty:

\begin{itemize}
\item SYN007.014 is perhaps inherently quite difficult because it involves many atoms with many equivalences.  Since $\leftrightarrow$ is treated as an abbreviation, the resulting formula is very large.
\item $\Pro$'s attacks are sometimes premature: for some formulas a solution can be found more quickly if an attack is delayed.
\end{itemize}

In light of the second observation, we considered an additional restriction on search:

\begin{description}
\item[Prefer-Defense] If $\Pro$ can defend, then he does defend.  (If multiple defenses are available, the choice is arbitrary.)
\end{description}

With this constraint, all the SYN problems of the ILTP library (except SYN007.014) are solvable, each in less than a minute (and most within several seconds).

Evaluating {\kuno} on properly first-order problems makes clear some difficulties with the naive depth-first approach currently implemented.  We can report, though, that working at the first-order level reveals challenges not revealed with the propositional SYN problems.  Namely, \textbf{Prefer-Defense} constraint cannot be rigidly applied; doing so leads to incompleteness.  A simple example (a modification of SYJ001+1.002) illustrates the difficulty:

\[
(\exists y \forall x (p(x) \wedge q(y))) \rightarrow (\forall x \exists y (p(x) \wedge q(y))).
\]

When playing a dialogue game for this formula, $\Opp$ begins by asserting the existential in the antecedent.  It is essential that $\Pro$ attacks this formula even though he could defend choose to defend against the initial attack by asserting the consequent.  The difficulty, intuitively, is that if $\Pro$ begins by defending, $\Opp$ can pin him down by attacking a formula that, from a constructive point of view, is weaker than what was initially given to $\Pro$.  To solve this, one apparently has to relax the constraint imposed by \textbf{Prefer-Defense}.

\section{Conclusion and future work}\label{sec:conclusion}

{\kuno} currently works in the equality-free fragment of intuitionistic logic.  This is, for many first-order theorem proving problems, a rather serious restriction that ought to be remedied.  At the moment, if a problem has equality in it at all, {\kuno} returns the SZS status \verb+Inappropriate+ to signal that it cannot deal with the problem.  A fairly easy remedy for such a gap would be to be preprocess (using, e.g., {\tptpfourx}) any problem involving equality so that appropriate equality axioms are present.  There seems to be no good treatment, from a dialogical point of view, of $\bottom$ and $\top$; {\kuno} rejects as inappropriate any TPTP problem that has these distinguished atomic formulas.  Still, one could presumably replace any occurrence of them by $p \rightarrow p$ and $p \wedge \neg p$, respectively.

By viewing dialogue games as an infrastructure for proof search, one quickly encounters the issue of redundancy.  Various strategies can in general be developed starting from an initial formula, but they can differ from one another in immaterial ways.  By contrast with resolution calculi, the notion of redundancy seems to be underdeveloped within the dialogue game framework.  Inspiration and ideas may come from other analytic search methods for intuitionistic logic, and perhaps even from game theory in general.


Intuitively, when E is present, $\Opp$ is more constrained than when E is absent.  One intuitively expects, then, that the E rules are favorable when the problem is to search for a winning strategy, that is, to determine intuitionistic validity.  When E is absent (that is, when playing according to the D rules), $\Opp$ has more options (at least sometimes, in general), so $\Pro$ intuitively faces a greater risk of losing ($\Opp$ might be able to pursue a more hostile ``line of reasoning'' against $\Pro$).  In this spirit, one could conduct further experiments that focus on non-theorems, rather than trying to verify theoremhood.  Such work could contribute to a better understanding of when dialogues go wrong~\cite{negri2014proofs}.

\bibliographystyle{plain}
\bibliography{paar}
\end{document}